%%%%%%%%%%%%%%%%THIS IS AN AMS-TEX DOCUMENT%%%%%%%%%%%%%%%%%%%%%%%%%%%%%%

\magnification=\magstep1
\vsize=22truecm
\input amstex
\documentstyle{amsppt}
\leftheadtext{E. Makai, Jr., J. Zem\'anek}
\rightheadtext{Nice connecting paths}
\topmatter
\title Nice connecting paths in connected components of sets of algebraic
elements in a Banach algebra\endtitle
\author 
{\centerline{Endre Makai, Jr., Budapest and Jaroslav Zem\'anek, Warszawa}}
%\endauthor
%\address 
\vskip.5cm
{\centerline{\it{Dedicated to the 90th anniversary of Professor Miroslav
Fiedler,}}}
{\centerline{\it{from two grateful participants in mathematical olympiads}}}
\vskip.5cm
{\centerline{MTA Alfr\'ed R\'enyi Institute of Mathematics,}}
{\centerline{H-1364 Budapest, Pf. 127, Hungary}}
{\centerline{{\rm{http://www.renyi.mta.hu/\~{}makai}}}}
\vskip.1cm
{\centerline{Institute of Mathematics, Polish Academy of Sciences,}}
{\centerline{00-656 Warsaw, \'Sniadeckich 8, Poland}}
%\endaddress
%\email 
\vskip.1cm
{\rm{makai.endre\@renyi.mta.hu, zemanek\@impan.pl}}
%\endemail
%\thanks *Research (partially) supported by Hungarian National Foundation for 
%Scientific Research, grant nos. ...\endthanks
\endauthor
%\keywords Banach algebras, $C^*$-algebras, (self-adjoint) idempotents,
%connected components of (self-adjoint) algebraic elements,
%(local) pathwise connectedness, similarities,
%analytic paths, polynomial paths, polygonal
%paths, centre of a Banach algebra, distance of connected components\endkeywords
%\subjclass {\it 2010 Mathematics Subject Classification.} Primary: 46H20, 
%Secondary: 46L05\endsubjclass
\abstract 
Generalizing earlier results about the set of idempotents in a Banach algebra,
or of self-adjoint idempotents in a $C^*$-algebra, we announce constructions of
nice connecting paths in
the connected components of the set of elements in a Banach algebra, or of
self-adjoint elements in a $C^*$-algebra, that satisfy a given polynomial
equation, without multiple roots. In particular, we will prove that in the
Banach algebra case every such non-central
element lies on a complex line, all of whose points satisfy the given equation.
We also formulate open questions.
\endabstract
\endtopmatter\document



\head 
1. Introduction
\endhead

\newpage

Let $A$ be a 
unital complex Banach algebra. Sometimes we will assume that moreover
$A$ is a $C^*$-algebra.

We let
$$
E(A):=\{ a \in A \mid a^2=a \}
$$
be the set of {\it{idempotents}} of $A$, and 
$$
S(A):= \{ a \in A \mid a^2=a=a^* \}
$$ 
be the set of {\it{self-adjoint idempotents}} for the $C^*$-algebra case.

The {\it{connected components of $E(A)$ and of $S(A)$}} 
have been investigated
by many authors. To some of them we will refer later at the respective
theorems. An ample literature is given in [AMMZ]. 

{\it{Let}} 
$$
p( \lambda ) := \prod _{i=1}^n (\lambda - \lambda _i)
$$
{\it{be a polynomial over ${\Bbb{C}}$, with all $\lambda _i$'s distinct. In the
$C^*$-algebra case, when considering self-adjoint elements, 
we will assume that all
$\lambda _i$'s are real.}} 
(In fact, if $q(\lambda ) := \prod \{ (\lambda - \lambda _i) \mid 1 \le i \le
n,\,\,\lambda _i \in
{\Bbb{R}} \} $, then $p(a)=0$ and $a=a^*$ imply $q(a)=0$. Thus below we could
use $q( \lambda )$ rather than $p(\lambda )$.)
The $\lambda _i$'s are fixed throughout this paper.

We write 
$$
E_p(A):= \{ a \in A \mid p(a)=0 \} ,
$$
and 
$$
S_p(A) := \{ a \in A \mid p(a)=0,\,\, a=a^* \}
$$
for the $C^*$-algebra case.
Then $E(A)$ and $S(A)$ are special cases of $E_p(A)$ and $S_p(A)$: namely, for
$p(\lambda ) := \lambda (\lambda -1)$.

We say that $\{ e_1, \ldots , e_n \} \subset A$ is a {\it{partition of
unity}}, or in the $C^*$-algebra case that 
$\{ e_1, \ldots , e_n \} \subset A$
is a {\it{self-adjoint partition of unity,}} if
$$
\cases
\{ e_1, \ldots , e_n \} \subset E(A), {\text{ or }}
\{ e_1, \ldots , e_n \} \subset S(A), \\
{\text{and }}
e_ie_j=0 {\text{ for }} 1 \le i,j \le n {\text{ and }} i \ne j, \\
{\text{and }} \sum _{i=1}^n e_i =1 .
\endcases
$$

\newpage

The detailed proofs of the statements announced 
in Section 2 will be published in [MZ]. The idea of this development originates
from personal conversations of the authors at the conference
Operator Theory and Applications: Perspectives and Challenges, held
in Jurata (Hel), Poland, March 18--28, 2010, 
and from the 2011 lecture by the first named author [Mak].

%%%%%%%%%%%%%%%%%%%%%%%%%%%%%%%%%%%%%%%%%%%%%%%%%%%%%%%%%%%%%%%%%%%%%%%%%%%%

\head
2. Theorems
\endhead

The ``only if'' part of the following Proposition 1 comes from the Riesz
decomposition theorem.

\proclaim{Proposition 1}
Let $A$ be a unital complex Banach algebra ($C^*$-algebra).
Let $a \in A$. Then $a \in E_p(A)$ ($a \in S_p(A)$)
if and only if there exists a (self-adjoint) partition of unity $\{ e_1,
\ldots , e_n \} $ such that
$$
a=\sum _{i=1}^n \lambda _i e_i .
$$
In the ``only if'' part, for $a \in E_p(A)$ (for 
$a \in S_p(A)$) one can choose the $e_i$'s as polynomials
of $a$, with complex (real) coefficients, which depend only on the $\lambda
_i$'s.
\endproclaim

This representation 
provides the tool for reducing questions about $E_p(A)$ (about $S_p(A)$)
to those about $E(A)$ (about $S(A)$). Of course, for the respective proofs for
$E_p(A)$ (for $S_p(A)$) one has to work still substantially. As an
illustration, we include a sketch of proof of Theorem 7 in Section 3.

The distinctness of the $\lambda _i$'s is essential in order that $a$ should
have such a simple form. For $T \in A:=B(l^2\oplus l^2)$, having a block
matrix form $(T_{ij})_{i,j=1}^2$, which is subdiagonal (i.e., strictly lower
triangular), we have $T^2=0$,
but $T_{21} \in B(l^2)$ can be as complicated as an element of $B(l^2)$ can be.

%%%%%%%%%%%%%%%%%%%%%%%%%%%%%%%%%%%%%%%%%%%%%%%%%%%%%%%%%%%%%%%%%%%%%%%%%%%

A {\it{path in a topological space}} $X$ is a continuous map $f: [0,1] \to X$.
We will say that $f(0),f(1) \in X$ {\it{are connected by this path}} $f$. 
By a small abuse of language we will also say that $f([0,1]) \subset X$ 
{\it{is a path in}} $X$ (e.g., for polygonal paths). A topological space
$X$ is {\it{pathwise
connected}} if any two of its points are connected by a path in $X$.
A topological space
$X$ is {\it{locally pathwise connected}} if each point $x \in X$ has a base of
(not necessarily open) neighbourhoods
consisting of pathwise connected sets.

%%%%%%%%%%%%%%%%%%%%%%%%%%%%%%%%%%%%%%%%%%%%%%%%%%%%%%%%%%%%%%%%%%%%%%%%%%%

\proclaim{Theorem 2}
Let $A$ be a unital complex Banach algebra and $C$ a connected component of
$E_p(A)$. Then $C$ is a relatively open subset of $E_p(A)$. Further, $C$ is
locally pathwise connected via each of the following types of paths:
\newline
1) similarity via an exponential function, i.e., $t \mapsto e^{-ct}ae^{ct}$;
\newline
2) a polynomial path of degree at most three;
\newline
3) a polygonal path of $n$ segments.
\endproclaim

\newpage

For $E(A)$, relative openness of $C$ was proved by J. Zem\'anek
[Ze],
1) was proved by J. Zem\'anek [Ze], 2) was proved by J. Esterle
[Es] and M. Tr\'emon [Tr85], 3) was proved by Z. V. Kovarik [Ko]
(cf. also [Ze]).

%%%%%%%%%%%%%%%%%%%%%%%%%%%%%%%%%%%%%%%%%%%%%%%%%%%%%%%%%%%%%%%%%%%%%%%%%%

\proclaim{Theorem 3}
Under the hypotheses of Theorem 2, $C$ is pathwise connected via each of the
following types of paths:
\newline
1) similarity via a finite product of exponential functions, i.e., 
$t \mapsto e^{-c_mt}\ldots e^{-c_1t}a$
\newline
$e^{c_1t} \ldots e^{c_mt}$;
\newline
2) a polynomial path;
\newline
3) a polygonal path.
\newline
In fact, there is a path satisfying 1) and 2) simultaneously.
\endproclaim

For $E(A)$, 1) was proved by J. Zem\'anek [Ze], 2) was proved by J. Esterle
[Es] and M. Tr\'emon [Tr85], 3) was proved by Z. V. Kovarik [Ko]
(cf. also [Ze]), and the last sentence was proved by [Es] and [Tr85].

%%%%%%%%%%%%%%%%%%%%%%%%%%%%%%%%%%%%%%%%%%%%%%%%%%%%%%%%%%%%%%%%%%%%%%%%%%

\vskip.1cm

{\bf{Problem.}} Does there exist a uniform bound on the ``minimum degree'' of
these polynomial connections, possibly depending on $n$, 
valid for all Banach algebras? Does such a bound exist, depending on $n$ and
on $A$ (or even on $C$)? 
Even the case of a uniform bound for polynomial connections of
idempotents is open, even if we allow dependence of the bound on $A$ (or even 
on $C$). 
For some particular cases, see [Tr85] and [MZ89]. ([Tr95] announced a further
partial result, but his proof seems to be incorrect.)

Even the ``simplest'' case $A:=B(l^2)$ is open. (The case $A=:B({\Bbb{C}}^n)$ is
solved positively by [Tr85], the uniform bound being $3$, which is sharp. 
Here the connected
components of $E(A)$ consist of the projections of the same rank.)
For $A=B(l^2)$, the connected components
of $E(A)$ are $\{ e \in A \mid {\text{dim}}\,N(e) = \alpha
,\,\,{\text{dim}}\,R(e) = \beta \} $, where $0 \le \alpha , \beta \le \aleph
_0$ are cardinalities with $\alpha + \beta = \aleph _0$, cf. [AMMZ] 
($N(\cdot )$ is the null-space and $R( \cdot )$ is the range). By [MZ89],
for $\min \{ \alpha , \beta \} < \aleph _0$, in the 
respective connected component there
exists an at most third degree polynomial path between any two elements of that
component. But even the case $\alpha = \beta = \aleph _0$ here is open.

\vskip.1cm

%%%%%%%%%%%%%%%%%%%%%%%%%%%%%%%%%%%%%%%%%%%%%%%%%%%%%%%%%%%%%%%%%%%%%%%%%%

\proclaim{Theorem 4}
Let $A$ be a unital complex $C^*$-algebra, and $C$ a connected component of
$S_p(A)$. Then $C$ is a  relatively open subset of $S_p(A)$. Further, $C$ is
locally pathwise connected by similarities via exponential functions, i.e.,
$t \mapsto e^{-ict} a e^{ict}$, where additionally $c=c^*$.
\endproclaim

For $S(A)$, Theorem 4 was proved by S. Maeda [Mae] (cf. also [Ze]). 

%%%%%%%%%%%%%%%%%%%%%%%%%%%%%%%%%%%%%%%%%%%%%%%%%%%%%%%%%%%%%%%%%%%%%%%%%%%%

\proclaim{Theorem 5}
Under the hypotheses of Theorem 4, $C$ is pathwise connected by similarities
via finite products of exponential functions, i.e., $t \mapsto e^{-ic_mt}
\ldots e^{-ic_1t} a e^{ic_1t} $
\newline
$\ldots e^{ic_mt}$, where
additionally $c_1=c_1^*, \ldots , c_m=c_m^*$.
\endproclaim

\newpage

For $S(A)$, Theorem 5 was proved by S. Maeda [Mae] (cf. also [Ze]). 

For the $C^*$-algebra case, the  analogues of 2) and 3) from Theorems 2 and 3
are false for $S_p(A)$. In fact, already the connected component of 
$S \left( B({\Bbb{C}}^2) \right )$ 
consisting of all rank-one orthogonal projections
does not contain any non-constant polynomial path.  
(The connected components of $S\left( B({\Bbb{C}}^n) \right) $ consist of the
orthogonal projections of the same rank.)

%%%%%%%%%%%%%%%%%%%%%%%%%%%%%%%%%%%%%%%%%%%%%%%%%%%%%%%%%%%%%%%%%%%%%%%%%%%%%

\proclaim{Theorem 6}
Let $A$ be a unital complex Banach algebra ($C^*$-algebra).
Let $a \in E_p(A)$ (let $a \in S_p(A)$). Then $a$ belongs to the centre of $A$ 
if and only if its connected component in $E_p(A)$ (in $S_p(A)$) is $\{ a \} $.
\endproclaim

Theorem 6 for $E(A)$ was proved by J. Zem\'anek [Ze], for $S(A)$ by S. Maeda
[Mae]. In Theorem 6, of course, the ``only if'' part for $S_p(A)$
follows from the ``only if'' part for $E_p(A)$.

%%%%%%%%%%%%%%%%%%%%%%%%%%%%%%%%%%%%%%%%%%%%%%%%%%%%%%%%%%%%%%%%%%%%%%%%%%%

\proclaim{Theorem 7}
Let $A$ be a unital complex Banach algebra, and $C$ a connected
component of $E_p(A)$. If $C$ is disjoint from the centre of $A$, then any
element of $C$ belongs to a complex line entirely contained in $C$.
In particular, $C$ is unbounded.
\endproclaim

For $E(A)$, Theorem 7 was proved by J. Zem\'anek [Ze].

In the $C^*$-algebra case even the entire $S_p(A)$ has a distance $\max \{ |
\lambda _i| \mid 1 \le i \le n \} $ from $0$, so the analogue of Theorem 7 
for $S_p(A)$ is false.

%%%%%%%%%%%%%%%%%%%%%%%%%%%%%%%%%%%%%%%%%%%%%%%%%%%%%%%%%%%%%%%%%%%%%%%%%%%

\vskip.1cm

Theorem 6 and Theorem 7 yield the next Corollary 8.

\proclaim{Corollary 8}
Let $A$ be a unital complex Banach algebra. Then $E_p(A)$ is a union
of its isolated points and of complex lines. $\blacksquare $
\endproclaim

%%%%%%%%%%%%%%%%%%%%%%%%%%%%%%%%%%%%%%%%%%%%%%%%%%%%%%%%%%%%%%%%%%%%%%%%%%%

\proclaim{Theorem 9}
There exists an explicit
constant $c(\lambda _1, \ldots , \lambda _n)>0$ (depending on
$\lambda _1, \ldots , \lambda _n \in {\Bbb{R}}$, and invariant under any map 
$(\lambda _1, \ldots , \lambda _n) \mapsto  (a+b\lambda _1, \ldots , 
a+b\lambda _n)$ with $a,b \in {\Bbb{R}}$ and $b \ne 0$) 
such that the following holds.  
If $A$ is a unital complex $C^*$-algebra, and $C_1,C_2$ are distinct connected
components of $S_p(A)$, then the distance of $C_1$ and $C_2$ is at least 
$c(\lambda _1, \ldots , \lambda _n) \cdot \min \{ | \lambda _i - \lambda _j|
\mid 1 \le i,j \le n,\,\, i \ne j \} $.
\endproclaim

%%%%%%%%%%%%%%%%%%%%%%%%%%%%%%%%%%%%%%%%%%%%%%%%%%%%%%%%%%%%%%%%%%%%%%%%%%%%%

{\bf{Conjecture.}}
Let $A$ be a unital complex Banach algebra ($C^*$-algebra) and $C_1,C_2$
distinct connected components of $E_p(A)$ (of $S_p(A)$).  
Then the distance of $C_1$ and $C_2$ is at least $\min \{ | \lambda _i -
\lambda _j| \mid 1 \le i,j \le n,\,\, i \ne j \} $.

\vskip.1cm

For $n=2$ this conjecture is equivalent to the statement 
that this distance for $E_p(A):=E(A)$ (for $S_p(A):=S(A)$) 
is at least $1$, which is due to
J. Zem\'anek [Ze] (due to  S. Maeda [Mae]). For $n \ge 3$ we do not even know
whether this distance for the Banach algebra case is positive.

\newpage

If true, this conjecture would be sharp, for any Banach algebra: 
consider $\lambda _i \cdot 1$ and $\lambda _j \cdot 1$.

The Conjecture for the case of $S_p(A)$ would follow from the Conjecture in
the case of $E_p(A)$. 
In fact, different connected components of $S_p(A)$ are subsets of 
different connected components of $E_p(A)$, by [BFML], Section 1,
Applications, 2), also taking into consideration our Proposition 1 and Theorem
3.

%%%%%%%%%%%%%%%%%%%%%%%%%%%%%%%%%%%%%%%%%%%%%%%%%%%%%%%%%%%%%%%%%%%%%%%%%%%%%%

\head 
3. A proof 
\endhead

\demo{Proof of Theorem 7 from Theorem 3 and Theorem 6}
If $C$ is disjoint from the centre, then by Theorem 6 it has
more than one elements. 
Let $a_0 \in C$ be an arbitrary element of $C$, 
and let $a_1 \in C$, with $a_1 \ne a_0$. Then, by Theorem 3, 3), 
there exists a non-constant polygonal path connecting
$a_0$ to $a_1$ in $C$. 
Its first non-constant segment (counted from $a_0$)
is the graph of a non-constant polynomial of degree $1$, say of
$$
\lambda \mapsto a_0+b \lambda , {\text{ from }} [0,1] {\text{ to }}
C\,\,(\subset E_p(A) \subset A) .
$$ 
Hence
$$
b \ne 0 {\text{ and we have for all }} \lambda \in [0,1] {\text{ identically }} 
p(a_0 + b \lambda )=0 .
\tag 1
$$ 
Then the equation in
\thetag{1} is a polynomial equation, with coefficients from $A$ and 
of degree at most $n$, 
for $\lambda \in {\Bbb{C}}$. (Attention: here the coefficient of $\lambda ^n$
is $b^n$, which may be $0$ even for $b \ne 0$.)

{\it{We make an indirect assumption. If the polynomial}} 
$$
{\Bbb{C}} \ni \lambda \mapsto p(a_0+b \lambda ) \in A
\tag 2
$$ 
{\it{were not identically $0$ for all}} $\lambda \in {\Bbb{C}}$, 
then for some $\lambda _0 \in {\Bbb{C}}$ we would have 
$$
p(a_0+b \lambda _0) \ne 0 .
$$
Then for some continuous linear functional $a'$ on $A$ we would have 
$$
\langle p(a_0+b \lambda _0), a' \rangle \ne 0 .
$$ 
The polynomial 
$$
{\Bbb{C}} \ni \lambda \mapsto 
\langle p(a_0+b \lambda ), a' \rangle \in {\Bbb{C}}
\tag 3
$$ 

\newpage

is a ${\Bbb{C}}$-valued polynomial on ${\Bbb{C}}$ of degree at most $n$, 
which would not vanish at $\lambda _0 \in {\Bbb{C}}$. Hence 
{\it{the polynomial in
\thetag{3} would have at most $n$ distinct roots}}. 

However, by \thetag{1} we have that {\it{the polynomial in \thetag{3}
vanishes for all $\lambda \in [0,1]$ identically}}. 
This is a 
contradiction, showing that our indirect assumption is false. 

That is, the polynomial in \thetag{2} is
identically $0$ for all $\lambda \in {\Bbb{C}}$. In other words, 
for all $\lambda \in {\Bbb{C}}$ we have 
$$
p(a_0 + b \lambda ) = 0, \,\,{\text{ i.e., }}\,\, a_0 + b \lambda  \in E_p(A),
$$
which implies by connectedness of ${\Bbb{C}}$ 
that for all $\lambda \in {\Bbb{C}}$ we have even
$$
a_0 + b \lambda  \in C .
$$

Since by \thetag{1} $b \ne 0$, we see that 
$$
C {\text{ contains a complex line passing through 
its arbitrary point }} a_0. 
$$
$\blacksquare $
\enddemo 

%%%%%%%%%%%%%%%%%%%%%%%%%%%%%%%%%%%%%%%%%%%%%%%%%%%%%%%%%%%%%%%%%%%%%%%%%%%%%%

{\bf{Acknowledgement.}} The authors are grateful to the organizers of several
conferences, where this material took its shape. In particular, to the
organizers of the
conference Operator Theory and Applications: Perspectives and Challenges, 
in Jurata (Hel), Poland, March 18--28, 2010, of the
6th Linear Algebra Workshop, in
Kranjska Gora, Slovenia, May 25--June 1, 2011, and of the Sz.-Nagy Centennial
Conference, in Szeged, Hungary, June 24--28, 2013.

%%%%%%%%%%%%%%%%%%%%%%%%%%%%%%%%%%%%%%%%%%%%%%%%%%%%%%%%%%%%%%%%%%%%%%%%%%%%%%

\end